\documentclass[12pt,a4paper]{article}

\usepackage[utf8]{inputenc}
\usepackage[english]{babel}

\usepackage{amsmath,amssymb,amsfonts,mathrsfs}

\usepackage[amsmath,thmmarks]{ntheorem}

\numberwithin{equation}{section}

\newtheorem{theorem}{Theorem}[section]
\newtheorem{corollary}[theorem]{Corollary}
\newtheorem{lemma}[theorem]{Lemma}
\newtheorem{proposition}[theorem]{Proposition}
\newtheorem{remark}[theorem]{Remark}
\newtheorem{definition}[theorem]{Definition}
\newtheorem{notation}[theorem]{Notation}

\theoremstyle{nonumberplain}
\theorembodyfont{\normalfont}
\theoremsymbol{\ensuremath{\square}}
\newtheorem{proof}{Proof}

\usepackage{hyperref}

\title{Loose ear decompositions and their applications to right-angled Artin groups}
\author{Max Gheorghiu}
\date{}

\usepackage[margin=2.5cm]{geometry}

\usepackage{tikz-cd}
\usetikzlibrary{decorations.pathmorphing}

\begin{document}

\setlength{\parindent}{0cm}

\maketitle

\begin{abstract}
We characterize planar graphs and graph minors among other graph theoretic notions in terms of right-angled Artin groups (RAAGs). For this, we determine all sets of elements in RAAGs with ears as underlying graphs that are exactly the sets of vertex generators. Generalizing ear decompositions of graphs to loose ear decompositions, we characterize both decompositions in terms of RAAGs. The desired results follow as applications of loose ear decompositions of RAAGs.
\end{abstract}

\tableofcontents

\section{Introduction}

We assume that every graph $\Gamma$ is finite, undirected and simple. Denote its vertex set by $V(\Gamma)$, its edge set by $E(\Gamma)$ and define its associated right-angled Artin group (RAAG) as the group
\[ A(\Gamma) := \langle v \in V(\Gamma) \mid [v, w] = 1 \text{ whenever } \lbrace v,  w \rbrace \in E(\Gamma) \rangle \, . \]
Adapting terminology from~\cite[p.~478]{kob22}, we call a generator arising from a vertex $v \in V(\Gamma)$ a vertex generator and a generating set as in the above group presentation a vertex basis. The graph $\Gamma$ is called the underlying graph of $A(\Gamma)$. The spirit of this paper is influenced by a proposed dictionary between graph theoretic terms and group theoretic terms that T.\! Koberda introduces in~\cite{kob22}. \\

The goal of this dictionary is to find for every property of a graph a corresponding algebraic property of a RAAG where no mention of generators or relations ought to be made for the latter. Geometrically, one can interpret this as a coordinate-free description. One can determine the underlying graph $\Gamma$ from the cohomology algebra of $A(\Gamma)$ as in~\cite[Theorem~15.2.6]{kob22} or by the elementary theory of $A(\Gamma)$ as in~\cite{cas21}. A functorial way is demonstrated in~\cite{gro23} which also provides an algorithm to determine $\Gamma$ from $A(\Gamma)$. However, the aim is to find non-trivial entries in this dictionary. To mention a few examples from the introduction of~\cite{kob22}, $\Gamma$ is $k$-colorable if and only if there is a surjection from $A(\Gamma)$ onto a $k$-fold direct product of free groups. Or $\Gamma$ possesses a non-trivial automorphism if and only if the outer automorphism group $\mathrm{Out}(A(\Gamma))$ contains a finite non-abelian group. We refer to~\cite{kob22} for more existing entries in this dictionary. \\

The graph theoretic notions we are mainly interested in are graph subdivisions, planar graphs and graph minors that we define as follows. An edge subdivision in a graph $\Gamma$ is defined by inserting a new vertex $v$ to $\Gamma$ and replacing an edge $\lbrace u, w \rbrace$ by the two edges $\lbrace u, v \rbrace$ and $\lbrace v, w \rbrace$~\cite[p.~305]{gro19}. A graph subdivision is defined as a finite sequence of edge subdivision~\cite[p.~306]{gro19}. A graph is planar if it can be drawn in the plane without edge-crossings~\cite[p.~298]{gro19}. Lastly, graph minors can be thought as a generalization of the notion of a subgraph. More rigorously, $X$ is a minor of $Y$ if there is $S \subseteq V(Y)$ and a function $f: S \rightarrow V(X)$ such that $f^{-1}(x) \leq Y$ is a connected subgraph and there is an edge between $f^{-1}(x)$ and $f^{-1}(x')$ if there is an edge between $x, x' \in V(X)$~\cite[p.~20]{die17}. Crucial for us are ears, namely graphs that are defined to consist of a single path or a single cycle~\cite[p.~252]{schr04}. To characterize the above graph theoretic notions in terms of RAAGs, we determine all sets of elements in RAAGs with ears as underlying graphs that are exactly the sets of vertex generators and call them pearl chains. \\

It is asked on page~513 in~\cite{kob22} whether one can characterize graph subdivisions in terms of RAAGs. Using pearl chains, we find for every edge subdivision a corresponding homomorphism between two RAAGs that we call a smoothing homomorphism. Hence we prove

\begin{lemma}
(= Lemma~\ref{lem:subdivisions}) Let $A(\Gamma)$ and $A(\Lambda)$ be RAAGs. Then $\Gamma$ is a graph subdivision of $\Lambda$ if and only if $A(\Gamma) \cong A(\Lambda)$ or there exists a finite sequence of smoothing homomorphisms $\lbrace \varphi_{i}: A(\Gamma_{i}) \rightarrow A(\Gamma_{i-1}) \rbrace_{i = 1}^m$ such that $A(\Gamma_m) = A(\Gamma)$ and $A(\Gamma_0) \cong A(\Lambda)$.
\end{lemma}

This lemma is used to settle Question~15.5.2 in~\cite{kob22} asking to characterize planar graphs in terms of RAAGs. This in turn requires to define loose ear decompositions of RAAGs in which pearls are elements that are part of pearl chains.

\begin{theorem}
(= Theorem~\ref{thm:planarity}) Let $A(\Gamma)$ be a RAAG. Assume we are given a sequence of smoothing homomorphisms $\lbrace \varphi_{i}: A(\Gamma_{i}) \rightarrow A(\Gamma_{i-1}) \rbrace_{i = 1}^m$ such that $A(\Gamma_m) = A(\Gamma)$ and such that there is no smoothing homomorphism from $A(\Gamma_0)$ to any other RAAG. Then the graph $\Gamma$ is planar if and only $A(\Gamma_0)$ does not contain $\mathbb{Z}^5$ as a subgroup and does not possess a loose ear decomposition with pearls $\lbrace p_i, q_i \rbrace_{i = 1}^3$ such that the subgroups $P := \langle p_1, p_2, p_3 \rangle$ and $Q := \langle q_1, q_2, q_3 \rangle$ are isomorphic to $F_3$ and $\langle P, Q \rangle \cong F_3 \oplus F_3$.
\end{theorem}

Moving forward, minors are an important object in graph theory. Namely, a class of graphs $\mathcal{P}$ is called minor-closed if any minor of a graph in $\mathcal{P}$ is again contained in $\mathcal{P}$. A graph $H$ is a forbidden minor for $\mathcal{P}$ if any graph in $\mathcal{P}$ does not contain $H$ as a minor~\cite[p.~369]{die17}. According to~\cite[p.~374]{die17}, ``theorems that characterize a property $\mathcal{P}$ by a set of forbidden minors are doubtless among the most attractive results in graph theory''. By pearl chains, certain relevant operations on graphs correspond to homomorphisms between RAAGs. We call a homomorphism associated with deleting a vertex a pearl deletion, one associated with deleting an edge an abelianizing homomorphism and one associated with identifying two vertices an identification homomorphism.

\begin{theorem}
(= Theorem~\ref{thm:minors}) Let $A(\Gamma)$ and $A(\Lambda)$ be RAAGs. Then $\Lambda$ is a minor of $\Gamma$ if and only if there exists a sequence of pearl deletions $\lbrace A(A_i) \rightarrow A(A_{i+1}) \rbrace_{i = 1}^a$, one of identification homomorphisms $\lbrace A(B_j) \rightarrow A(B_{j+1}) \rbrace_{j = 1}^b$ and one of abelianizing homomorphisms $\lbrace A(C_k) \rightarrow A(C_{k+1}) \rbrace_{k = 1}^c$ such that $A(\Gamma) \cong A(A_1)$, $A(A_a) \cong A(B_1)$, $A(B_b) \cong A(C_c)$ and $A(C_1) \cong A(\Lambda)$.
\end{theorem}

Our characterizations pertain to the realm of T.\! Koberda's dictionary because they do not mention generators and relations explicitly. But they do not fit properly into the dictionary because they use generators and relations implicitly. More explicitly, we characterize vertex bases of ear RAAGs without stating upfront that they are generators, call these pearl chains and obtain our characterizations through them. Moreover, given a RAAG $A(\Gamma)$, it is unlikely that one could implement our characterizations into an algorithm to determine whether $\Gamma$ possesses any of these graph theoretic properties. In the meantime, in the paper~\cite{flo23} a characterization of $\Gamma$ in terms of $A(\Gamma)$ has been found through group cohomology such that it does not depend on any choice of generators or relations of the RAAG. Moreover, this characterization can be algorithmically implemented (which the authors of that paper call effective). In particular, this yields an algorithm determining from any finite presentation of $A(\Gamma)$ whether $\Gamma$ is planar. We wonder whether there is an (algorithmic) characterization of planarity in terms of RAAGs that does not resort to the underlying graph, but only uses algebraic means such as taking subgroups, quotients, elements, centralizers etc. \\

Our characterizations all depend on a generalization of ear decompositions of graphs that we call loose ear decompositions. We assume that every ear, path and cycle contains more than one vertex. In order to take over a term for paths to ears, we call one distinguished vertex of a cycle its end point and all other vertices inner vertices (thus, one can think of a cycle as path whose endpoints are identified). An ear decomposition of $\Gamma$ is a finite sequence of ears $\lbrace \Gamma_i \rbrace_{i = 1}^m$ contained in the graph $\Gamma$ and satisfying the following properties. The ear $\Gamma_1$ is a cycle. For any $i > 1$ only the endpoints of the ear $\Gamma_i$ are contained in the previous ears $\bigcup_{j = 1}^{i-1} \Gamma_j$, but none of its inner vertices (if there are any). Finally, $\Gamma = \bigcup_{i = 1}^m \Gamma_i$~\cite[p.~235--237]{gro19}, \cite[Definition~1]{schm16}. We define a loose ear decompositions of $\Gamma$ as a finite sequence of ears $\lbrace \Gamma_i \rbrace_{i = 1}^m$ in the same way, but without imposing that $\Gamma_1$ is a cycle or that the end points of $\Gamma_i$ have to be contained in $\bigcup_{j = 1}^{i-1} \Gamma_j$. \\

In the initial sections, we consider ears and call RAAGs with ears as underlying graphs ear RAAGs. In the second section, we characterize elements in almost all ear RAAGs that arise from degree $2$ vertices without specifying that these are generators (Lemma~\ref{lem:degtwovertices}). We determine in the third section when two such elements arise from adjacent vertices (Proposition~\ref{prop:adjacency}) as well as all elements that arise from degree $1$ vertices (Lemma~\ref{lem:degonevertices}). Subsequently, we characterize in the fourth section vertex bases in ear RAAGs that do not meet the conditions of the previous sections (Proposition~\ref{prop:specialcases}). We call these elements arising from degree $1$ or $2$ vertices pearls. In particular, pearls are vertex generators. In the fifth section, we characterize the sets of pearls in ear RAAGs that correspond exactly to vertex bases (Lemma~\ref{lem:pathpearlchains} and Lemma~\ref{lem:cyclepearlchains}) and call such sets pearl chains. This is non-trivial as there is no guarantee to obtain a vertex basis just by picking vertex generators. Although pearl chains do not fit into T.\! Koberda's dictionary, they still pertain to its realm because we do not state that they are actually generating sets. Using them, we also characterize ear decompositions and loose ear decomposition of a graph $\Gamma$ in terms of its associated RAAG $A(\Gamma)$ in this section. To be more explicit, a loose ear decomposition of $\Gamma$ is equivalent to a specific iterated amalgamated product of ear RAAGs that is isomorphic to $A(\Gamma)$ (Theorem~\ref{thm:looseears}). \\

In the sixth section, we showcase applications of loose ear decompositions and ear decompositions to RAAGs. First, we characterize graphs subdivisions and planar graphs in term of their associated RAAGs (Lemma~\ref{lem:subdivisions} and Theorem~\ref{thm:planarity}). Then we describe graph minors through RAAGs (Theorem~\ref{thm:minors}). This in turn allows us to determine other graph theoretic properties algebraically by their associated RAAGs. For instance, if $\Gamma$ is a graph, we provide conditions on the associated RAAG $A(\Gamma)$ determining when $\Gamma$ is a forest or an outerplanar graph (Proposition~\ref{prop:forests} and Lemma~\ref{lem:outerplanar}). Graph minors could be also used to obtain conditions on $A(\Gamma)$ for when $\Gamma$ is an apex graph and embeds into a given surface without edge-crossings. In our last applications, we restrict our attention to ear decompositions. Doing so, we find algebraic criteria for $A(\Gamma)$ determining when $\Gamma$ is $2$-edge-connected (Theorem~\ref{thm:twoedgeconn}), $2$-vertex-connected (Lemma~\ref{lem:twovertexconn}), a series-parallel graph (Lemma~\ref{lem:seriesparallel}), a factor-critical graph (Lemma~\ref{lem:factorcritical}) and what the maximal size of a join in it is (Proposition~\ref{prop:maxjoinsize}).

\section{Elements in ear RAAGs arsing from degree 2 vertices}

In this section, we characterize elements in ear RAAGs that arise from degree $2$ vertices without specifying that these are generators. For this we first need to characterize ear RAAGs among all RAAGs which necessitates the following piece of notation.

\begin{notation}
For a finitely generated group $G$, let $d(G)$ denote the minimal number of generators of $G$.
\end{notation}

\begin{proposition}
If a RAAG $A(\Gamma)$ does not split as non-trivial free product and for any $x \in A(\Gamma)$ we have $d(C(x)) \leq 3$, then $A(\Gamma)$ is an ear RAAG.
\end{proposition}

\begin{proof}
According to Proposition~15.3.6 in~\cite{kob22}, the degree of every vertex in $\Gamma$ is at most $2$. By Theorem~15.3.2 in~\cite{kob22}, the graph $\Gamma$ is connected. If $\Gamma$ has a vertex of degree $1$, then it is a path. If not, then it is a cycle.
\end{proof}

In the proofs of this paper, we also need an elementary result for which we have not found an explicit proof in the literature.

\begin{proposition}\label{prop:degsandvertices}
For any RAAG $A(\Gamma)$ we have $d(A(\Gamma)) = |V(\Gamma)|$.
\end{proposition}

\begin{proof}
Because we can define $A(\Gamma)$ by a group presentation using the graph $\Gamma$, we see that $|V(\Gamma)| \geq d(A(\Gamma))$. By this presentation, we observe for the abelianization that $A(\Gamma)_{ab} = \mathbb{Z}^{|V(\Gamma)|}$. Using the projection homomorphism $A(\Gamma) \rightarrow A(\Gamma)_{ab}$ we infer that $d(A(\Gamma)) \geq d(\mathbb{Z}^{|V(\Gamma)|}) = |V(\Gamma)|$.
\end{proof}

\begin{lemma}
If $A(\Gamma)$ is an ear RAAG with $\Gamma$ not isomorphic to the cycle on $4$ vertices $C_4$ and $x \in A(\Gamma)$ a nontrivial element with a non-abelian centralizer, then, up to conjugation,
\begin{itemize}
\item $x$ arises as a power of a degree $2$ vertex $v$ in $\Gamma$ and
\item $C(x) \cong \langle v \rangle \oplus F_2$.
\end{itemize}
\end{lemma}

\begin{proof}
Denote by $v_1, \dots ,v_n$ the vertices of $\Gamma$ where $v_i$ is adjacent to $v_{i+1}$. If $\Gamma$ is a path, we assume that $v_1$, $v_n$ are the endpoints and if $\Gamma$ is a cycle, then $v_1$ is adjacent to $v_n$ By abuse of notation, we shall denote the corresponding vertex basis of $A(\Gamma)$ also by $v_1, \dots, v_n$. Writing $x$ as a word in these generators, we would like to apply the centralizer theorem on page~44 from Servatius's paper~\cite{ser89} to $C(x)$. Following~\cite[p.~489]{kob22}, we call a word in the vertex basis $v_1, \dots, v_n$ reduced if it cannot be shortened using free reductions or commutations. Moreover a word is called cyclically reduced if it is remains reduced after any cyclic permutation of its letters. Since the statement of the lemma will not change under conjugation, we may assume that $x$ is cyclically reduced. Define the support of $x$ to be the subset of our vertex basis required to write $x$ as a reduced word and denote it by $\mathrm{supp}(x)$~\cite[p.~489]{kob22}. This notion is well defined in the sense that for reduced words $w_1 = w_2$ in $A(\Gamma)$ we have $\mathrm{supp}(w_1) = \mathrm{supp}(w_2)$. Denote by $\Lambda_1, \dots, \Lambda_k \leq \Gamma$ the subgraphs such that the support of $x$ can be written as a join $\mathrm{supp}(x) = \Lambda_1 \ast {} \dots {} \ast \Lambda_k$. For $i \in \lbrace 1, \dots, k \rbrace$ choose $u_i \in A(\Lambda_i)$ such that $x = \prod_{i = 1}^k u_i^{m_i}$ with $m_i$ maximal. By Servatius's centralizer theorem  we can write
\[ C(x) = \prod_{i = 1}^k \langle u_i \rangle \oplus \langle \mathrm{link}(\mathrm{supp}(x)) \rangle \, . \]
Because $C(x)$ is non-abelian, we see that $\mathrm{link}(\mathrm{supp}(x)) \neq \varnothing$. Moreover, if $P_2$ denotes the path on $2$ vertices, then $A(P_2) = \mathbb{Z}^2$ and $A(C_3) = \mathbb{Z}^3$ whence $\Gamma \notin \lbrace P_2, C_3 \rbrace$. Note that
\[ \mathrm{link}(v_i) \cap \mathrm{link}(v_{i+1}) = \varnothing \text{ and } \mathrm{link}(v_i) \cap \mathrm{link}(v_{i+2}) = \lbrace v_{i+1} \rbrace \]
because $\Gamma \neq C_4$. If $\Gamma$ is a cycle, we may take the index $i$ modulo $|V(\Gamma)|$. Hence $\mathrm{link}(\mathrm{supp}(x)) \subseteq \lbrace v_i, v_{i+1}, v_{i+2} \rbrace$. This leaves us with the following cases to consider:
\begin{itemize}
    \item Case 1: $x \in \langle v_i \rangle$ for a vertex $v_i$ that is not an endpoint.
    \item Case 2: $x \in \langle v_i \rangle$ for $v_i$ an endpoint if $\Gamma$ is a path.
    \item Case 3: $\mathrm{supp}(x) = \lbrace v_i, v_{i+2} \rbrace$.
\end{itemize}
Going through these cases, we notice that only Case 1 yields $C(x) \cong \langle v_i \rangle \oplus F_2$ while all other cases yield $C(x) \cong \mathbb{Z}^2$. We conclude that $x = v_i^m$ for $m \in \mathbb{Z} \setminus \lbrace 0 \rbrace$.
\end{proof}

\begin{lemma}\label{lem:degtwovertices}
Let $A(\Gamma)$ be an ear RAAG with $\Gamma \neq C_4$. Let $x \in A(\Gamma)$ be a nontrivial element with non-abelian centralizer such that $x$ generates the center of its centralizer. Then, up to conjugacy and inversion, $x$ arises from degree $2$ vertex in $\Gamma$.
\end{lemma}

\begin{proof}
By the previous lemma, there is up to conjugacy a degree $2$ vertex $v$ and an integer $m \in \mathbb{Z} \setminus \lbrace 0 \rbrace$ such that $x = v^m$. Note that any element in the centralizer of $C(x) = \langle v \rangle \oplus F_2$ is of the form $v^a$ for $a \in \mathbb{Z}$. Thus, $x = v^{\pm 1}$.
\end{proof}

\begin{corollary}
Let $x \in A(\Gamma)$ satisfy the condition of the Lemma~\ref{lem:degtwovertices}. If $A(\Gamma)/ \langle\! \langle x \rangle\! \rangle$ splits as a non-trivial free product, then $\Gamma$ is a path. Otherwise it is a cycle.
\end{corollary}

\begin{proof}
This follows from Theorem~15.3.2 in~\cite{kob22} and from the fact the the element $x$ corresponds to a degree $2$ vertex in $\Gamma$.
\end{proof}

\begin{definition}
A RAAG with a path (resp.\! cycle) as underlying graph is called a path RAAG (resp.\! cycle RAAG).
\end{definition}

\section{Elements in ear RAAGs arsing from degree 1 vertices}

In this section, we characterize all elements in path RAAGs that arise from the degree $1$ vertices as well as all pairs of elements in an ear RAAG that arise from pairs of adjacent degree $2$ vertices without specifying that these are generators. \\

For this, we set up a few preliminary definitions for ear RAAGs in accordance with the terminology in graph theory (where we replace the word ``vertex'' by ``pearl''). Namely, inner vertices and endpoints are fundamental notions for loose ear decompositions of graphs. Because loose ear decompositions of RAAGs are our main technical tool, we aim to introduce the notions of inner pearls and end pearls. If $A(\Gamma)$ is an ear RAAG, then the elements of $A(\Gamma)$ satisfying the conditions of the Lemma~\ref{lem:degtwovertices} are an example of what we call pearls. On the one hand, if $A(\Gamma)$ is a path RAAG, then these elements are also an example of inner pearls. On the other hand, if $A(\Gamma)$ is a cycle RAAG, then we call a distinguished pearl of $A(\Gamma)$ the end pearl. All other pearls in a cycle RAAG are called inner pearls. At the end of the next section we shall be able to define the above notions rigorously. For the time being, an element in $A(\Gamma)$ is called an inner pearl if it satisfies the conditions in Lemma~\ref{lem:degtwovertices}. The distinction between inner pearls and end pearls is important in the characterization of graph subdivisions and planarity through RAAGs (Lemma~\ref{lem:subdivisions} and Theorem~\ref{thm:planarity}) as well as in all applications of ear decompositions of RAAGs at the end of the sixth section.

\begin{proposition}\label{prop:adjacency}
Let $A(\Gamma)$ be an ear RAAG with $\Gamma \neq C_4$ and $x, y \in A(\Gamma)$ be inner pearls. Suppose that $[x, y] = 1$ and that there is $z \in A(\Gamma)$ such that  $[x, z] = 1$ and $[y, z] \neq 1$. Then $x$ and $y$ form part of vertex basis such that they arise from adjacent vertices.
\end{proposition}

\begin{proof}
Although $x$ and $y$ are vertex generators arising from degree $2$ vertices, they might not form part of any vertex basis together. Following Lemma~\ref{lem:degtwovertices}, we may identify $x$ with a vertex $v_i$ and write $y = av_j a^{-1}$ for an element $a \in A(\Gamma)$. Assume by contradiction that $a \notin C(x)$. Writing $a$ as a word in the vertex generators $\lbrace v_1, \dots, v_n \rbrace$, we may assume that it is reduced. If none of the last letters of $a$ is in $C(v_j)$ even after commutation, then $y = av_j a^{-1}$ is also a reduced word. Whether $x$ and $y$ form part of a vertex basis is invariant under conjugation. Since $x$ is invariant under conjugation by elements in $C(x)$ and we obtain $y$ by conjugating $v_j$ by $a$, we may further assume that none of the initial letters of the word $a$ is in $C(x)$ even after commutation. This implies that the word $[x, y] = v_i av_j a^{-1}v_i^{-1}av_j^{-1}a^{-1}$ is reduced. On page~489 of~\cite{kob22}, it is mentioned that the word problem in RAAGs can be solved by free reductions and commutations. This yields the contradiction that $[x, y] \neq 1$. Thus, conjugation by $a \in C(x)$ maps $\lbrace x, y \rbrace$ to the vertex generators $\lbrace v_i, v_j \rbrace$. Using that $[x, y] = 1$, $[x, z] = 1$ and $[y, z] \neq 1$, we conclude that $\Gamma$ is not isomorphic to the path on $3$ vertices $P_3$, that $j \in \lbrace i-1, i+1 \rbrace$ and that $x$ and $y$ arise indeed from adjacent degree $2$ vertices.
\end{proof}

Pearls such as in the above proposition are an example of what we call adjacent pearls. This notion is used to characterize all elements in path RAAGs that arise from degree $1$ vertices.

\begin{lemma}\label{lem:degonevertices}
Let $A(\Gamma)$ be an ear RAAG with underlying graph $\Gamma \neq C_4$, $x \in A(\Gamma)$ an element and $y, z \in A(\Gamma)$ adjacent inner pearls such that $C(y) = \langle x, y, z \rangle$ and such that $A(\Gamma)/ \langle\! \langle y \rangle\! \rangle$ splits as a free product with one factor being $\mathbb{Z}$. Then there is a vertex basis of  $A(\Gamma)$ containing $x$, $y$ and $z$ such that $x$ arises from a degree $1$ vertex in $\Gamma$, $y$ from its adjacent vertex and $z$ from the degree $2$ vertex adjacent to $y$.
\end{lemma}

\begin{proof}
Note that all the conditions on the elements $x$, $y$ and $z$ are invariant under conjugation. Since $y$ and $z$ are adjacent inner pearls, we may assume that $y$ arises from a degree $2$ vertex $v_i$ and $z$ from either $v_{i-1}$ or $v_{i+1}$. As $A(\Gamma)/ \langle\! \langle y \rangle\! \rangle$ splits as a non-trivial free product, $A(\Gamma)$ is a path RAAG. More specifically, there exist two paths $\Gamma'$ and $\Gamma''$ such that $\Gamma \setminus \lbrace y \rbrace = \Gamma' \sqcup \Gamma''$ and $A(\Gamma)/ \langle\! \langle y \rangle\! \rangle = A(\Gamma') \ast A(\Gamma'')$. We already know that there is also a splitting of the form $A(\Gamma)/ \langle\! \langle y \rangle\! \rangle = \mathbb{Z} \ast A$. By to the Grushko Decomposition Theorem~\cite[p.~168--169]{sta77}, we infer that $A(\Gamma') \cong \mathbb{Z}$ and $A(\Gamma'') \cong A$. Since isomorphic RAAGs have the same underlying graph~\cite[Theorem~15.2.6]{kob22}, the graph $\Gamma'$ consists of a single vertex. If we assume that $v_1$ and $v_n$ are the endpoints of $\Gamma$, then the only choices of $y$ inducing an isolated vertex in $\Gamma \setminus \lbrace y \rbrace$ are $y = v_2$ and $y = v_{n-1}$. Let us assume $y = v_2$. As $z$ is an inner pearl adjacent to $y$, we have $z = v_3$. \\

Up to the end of the proof, we restrict our attention to the subgroup $C(y)$. Since both $x \langle\! \langle y, z \rangle\! \rangle$ and $v_1 \langle\! \langle y, z \rangle\! \rangle$ generate $C(y)/ \langle\! \langle y, z \rangle\! \rangle \cong \mathbb{Z}$, there is $c = \pm 1$ with the property that $x \langle\! \langle y, z \rangle\! \rangle = v_1^c \langle\! \langle y, z \rangle\! \rangle$. In particular, the generator $x \langle\! \langle y, [v_1, v_3] \rangle\! \rangle$ of $C(y)/ \langle\! \langle y, [v_1, v_3] \rangle\! \rangle \cong \mathbb{Z}^2$ is of the form $v_1^c v_3^d \langle\! \langle y, [v_1, v_3] \rangle\! \rangle$. Note that the projection homomorphism
\[ C(y)/ \langle y \rangle \cong F_2 \rightarrow C(y)/ \langle\! \langle y, [v_1, v_3] \rangle\! \rangle \cong \mathbb{Z}^2 \]
is the homomorphism to the abelianization. According to Corollary~3.3 in~\cite{osb81}, there is, up to conjugacy, a unique vertex generator in $F_2$ mapping down to $x \langle\! \langle y, [v_1, v_3] \rangle\! \rangle$. For any cyclic permutation $f$ of the word $v_1^c v_3^d$ the element $f \langle y \rangle$ has this property. There is such $f$ together with a reduced word $a$ in the letters $\lbrace v_1^{\pm 1}, v_3^{\pm 1} \rbrace$ such that $afa^{-1} \langle y \rangle$ is a reduced word in $\lbrace v_1^{\pm 1} \langle y \rangle v_3^{\pm 1} \langle y \rangle \rbrace$ and $x \langle y \rangle = afa^{-1} \langle y \rangle$. Whether $\lbrace x \langle y \rangle, z \langle y \rangle \rbrace$ generates $C(y)/ \langle y \rangle$ is invariant under conjugation by $v_3^{\pm 1} \langle y \rangle$. Thus, we may assume that the initial letter of $a$ is not $v_3 \langle y \rangle$ or $v_3^{-1} \langle y \rangle$. If the word $a$ is nontrivial after this last conjugation, then we obtain the contradiction that $x \langle y \rangle$ and $z \langle y \rangle$ do not generate $C(y) / \langle y \rangle$. We conclude that there are integers $\gamma, \delta \in \mathbb{Z}$ such that $x \langle y  \rangle = v_3^{\gamma}v_1^c v_3^{\delta} \langle y \rangle$. In particular, there is $\eta \in \mathbb{Z}$ such that $x = v_3^{\gamma}v_1^c v_3^{\delta}v_2^{\eta}$. According to the treatment of automorphisms of RAAGs on page~498 of~\cite{kob22}, we may apply so-called dominated transvections mapping $v_1$ to $v_1 v_2^{\pm 1}$, $v_1 v_3^{\pm 1}$ or $v_3^{\pm 1}v_1$ and leaving all other vertex generators unchanged. Therefore, we may assume that $x = v_1^c$ for $c = \pm 1$, $y = v_2$ and $z = v_3$. Invoking again~\cite[p.~498]{kob22}, we may apply the so-called vertex inversion mapping $v_1$ to $v_1^{-1}$ and leaving all the other vertex generators invariant. Hence, $x = v_1$, $y = v_2$ and $z = v_3$.
\end{proof}

\section{Special cases of ear RAAGs}

We characterize in this section vertex bases in special cases of RAAGs without saying upfront that the elements are generators. This allows us at the end of the section to define pearls, inner pearls and end pearls rigorously. Lemma~\ref{lem:degtwovertices} only characterizes elements in an ear RAAG $A(\Gamma)$ arising from a degree $2$ vertex whenever $\Gamma$ is not isomorphic to the cycle on $4$ vertices $C_4$. Due to the proof of Proposition~\ref{prop:adjacency}, Lemma~\ref{lem:degonevertices} only characterizes elements in an ear RAAG $A(\Gamma)$ arising from a degree $1$ vertex whenever $\Gamma$ is not isomorphic to the path on $3$ vertices $P_3$. None of the above results characterize elements in the ear RAAGs $A(P_2) = \mathbb{Z}^2$ and $A(C_3) = \mathbb{Z}^3$ arising from the vertices of their underlying graphs.

\begin{proposition}\label{prop:specialcases}
\begin{enumerate}
    \item Two elements $x, y \in A(P_2) = \mathbb{Z}^2$ form a vertex basis if and only if $x$ and $y \langle x \rangle \in A(P_2)/ \langle x \rangle$ are not multiples of other elements.
    \item A vertex basis of $A(C_3) = \mathbb{Z}^3$ is formed by $\xi, \eta, \zeta$ if and only if $A(C_3)/ \langle \zeta \rangle \cong \mathbb{Z}^2$ and $\xi \langle \zeta \rangle, \eta \langle \zeta \rangle \in A(C_3)/ \langle \xi \rangle$ satisfy the same conditions as $x$ and $y$ above.
    \item A vertex basis of $A(P_3) = \mathbb{Z} \oplus F_2$ is formed by $\chi, a, b$ with $\chi$ corresponding to the inner vertex and $a$ and $b$ to the endpoints if and only if $\chi$ generates the center of its non-abelian centralizer $C(\chi)$ and $C(\chi) = \langle \chi, a, b \rangle$.
    \item Let $u_1, u_2, u_3, u_4$ be elements in the RAAG $A(C_4) = F_2 \oplus F_2$. Then they form a vertex basis if and only if for any index $i$ modulo $4$ the centralizer $C(u_i)$ is non-abelian, $C(u_i) = \langle u_{i-1}, u_i, u_{i+1} \rangle$ and $A(C_4)/ \langle\! \langle u_i \rangle\! \rangle$ is a RAAG with the property that $d(A(C_4)) = d(A(C_4)/ \langle\! \langle u_i \rangle\! \rangle)+1$.
\end{enumerate}
\end{proposition}

\begin{proof}
\textbf{(1)} If $x$ is a vertex generator, then $\mathbb{Z}^2/ \langle x \rangle \cong \mathbb{Z}$ is torsion-free, hence $x$ is not a multiple of a nontrivial element. Further, an element in $\mathbb{Z}$ is a vertex generator if and only if it is no multiple of another element. It remains to prove that $x$ is a vertex generator whenever it is no multiple of another element. If $s, t \in \mathbb{Z}^2$ form a generating set, then there are coprime integers $\alpha, \beta$ such that $x = {\alpha}s + {\beta}t$. According to Bézout's Theorem, there are integers $c, d$ such that ${\alpha}c+{\beta}d = 1$. Observing that $s = c({\alpha}s + {\beta}t)+{\beta}(ds - ct)$ and $t = d({\alpha}s + {\beta}t)-{\alpha}(ds - ct)$, we conclude that $x$ is a vertex generator. \\

\textbf{(2)} This follows from the first assertion. \\

\textbf{(3)} If $\chi, a, b$ form a vertex basis, then it is straightforward to show that they satisfy the other assertion. Let us prove the reverse implication. If $\chi$ was trivial, then $A(P_3) = C(\chi)$ was generated by $a$ and $b$ although $d(A(P_3)) = 3$ by Proposition~\ref{prop:degsandvertices}. Hence, $\chi$ is nontrivial. If $V(P_3) = \lbrace v_1, v_2, v_3 \rbrace$ with $v_2$ being the inner vertex, then we may assume by Lemma~\ref{lem:degtwovertices} that $\chi = v_2^{\pm 1}$. Because $A(P_3) = C(\chi) = \langle \chi, a, b \rangle$, we conclude that $\lbrace \chi, a, b \rbrace$ forms a vertex basis of $A(P_3)$. \\

\textbf{(4)} If $u_1, u_2, u_3, u_4$ form a vertex basis, then the other assertion follows from Proposition~\ref{prop:degsandvertices} and the centralizer theorem stated in the proof of Lemma~\ref{lem:degtwovertices}. Let us demonstrate the reverse implication. As before, we write the vertices of $C_4$ as $\lbrace v_1, v_2, v_3, v_4 \rbrace$. Up to conjugation, we may assume that $u_1$ is cyclically reduced. By the centralizer theorem and the fact that $u_1$ has a non-abelian centralizer, we infer that $\mathrm{link}(\mathrm{supp}(u_1)) \neq \varnothing$ and that we only need to consider the two cases
\begin{itemize}
    \item Case 1: $u_1 \in \langle v_i \rangle$.
    \item Case 2: $u_1 \in \langle v_i, v_{i+2} \rangle$.
\end{itemize}
Note that each $F_2$-factor of $A(P_4) = F_2 \oplus F_2$ is generated by a tuple of the form $\lbrace v_i, v_{i+2} \rbrace$. We may thus write $u_1 = (U_1, 1)$. If we only consider the subgroup $\langle v_i, v_{i+2} \rangle \cong F_2$, then its abelianization is given by $\langle v_i, v_{i+2} \rangle/ \langle\! \langle [v_i, v_{i+2}] \rangle\! \rangle \cong \mathbb{Z}^2$. According to Theorem~1.2 and Corollary~3.3 in R.\! P.\! Osborne and H.\! Zieschang's paper~\cite{osb81}, for every vertex generator in the abelianization $\mathbb{Z}^2$ there is a vertex generator in $F_2$ mapping down to it and all vertex generators in $F_2$ are of this form. If we assume by contradiction that $U_1$ is not a vertex generator, then either $U_1 \langle\! \langle [v_i, v_{i+2}] \rangle\! \rangle = 0$ or there are $0 \neq U \in \langle v_i, v_{i+2} \rangle/ \langle\! \langle [v_i, v_{i+2}] \rangle\! \rangle$ and $k \geq 2$ such that $U^k = U_1 \langle\! \langle [v_i, v_{i+2}] \rangle\! \rangle$. Let us suppose first that $U_1 \in \langle\! \langle [v_i, v_{i+2}] \rangle\! \rangle$. Then the projection homomorphism to the abelianization factors as $A(C_4) \rightarrow A(C_4)/ \langle\! \langle u_1 \rangle\! \rangle \rightarrow \mathbb{Z}^4$. However, this yields the contradiction
\[ 3 = d(A(C_4))-1 = d(A(C_4)/ \langle\! \langle u_1 \rangle\! \rangle) \geq d(\mathbb{Z}^4) = 4 \, . \]
Therefore $U^k = U_1 \langle\! \langle [v_i, v_{i+2}] \rangle\! \rangle$. Then the abelianization
\[ (A(C_4)/ \langle\! \langle u_1 \rangle\! \rangle)_{ab} \cong \big((F_2/ \langle\! \langle U_1 \rangle\! \rangle) \oplus F_2 \big)_{ab} \cong (\mathbb{Z}^2/ \langle U_1 \langle\! \langle [v_i, v_{i+2}] \rangle\! \rangle \, \rangle) \oplus \mathbb{Z}^2 \]
contains the torsion element $(U, 1)$. Since abelianizations of RAAGs are torsion-free, our assumption has led to a contradiction and $U_1 \in \langle v_i, v_{i+2} \rangle$ is a vertex generator. If $U_1$ together with $V$ generate the subgroup $\langle v_i, v_{i+2} \rangle \cong F_2$, then $u_1, (V, 1), v_{i+1}, v_{i+3}$ form another vertex basis of $A(C_4)$. If we write $u_2 = (W, U_2)$, then $u_2$ is contained in $C(u_1) = \langle u_1, v_{i+1}, v_{i+3} \rangle$ and satisfies the above two cases as $u_1$ does. Therefore, $u_2 = u_1^{\pm 1}$ or $u_2$ is a vertex generator in $\langle v_{i+1}, v_{i+3} \rangle \cong F_2$. Because $u_1$ arises from a degree $2$ vertex, $C(u_1) \cong A(P_3)$ and thus $d(C(u_1)) = 3$. Since $C(u_1) = \langle u_4, u_1, u_2 \rangle$, we deduce that $u_2$ is a vertex generator in $\langle v_{i+1}, v_{i+3} \rangle$. Note that the same holds true for the element $u_4$. Given that the tuple $\lbrace u_4 \langle u_1 \rangle, u_2 \langle u_1 \rangle \rbrace$ is a generating set of $\langle v_{i+1} \langle u_1 \rangle, v_{i+3} \langle u_1 \rangle \rangle = C(u_1)/ \langle u_1 \rangle$, we conclude that $\lbrace u_1, u_2, (V, 1), u_4 \rbrace$ forms another vertex basis of $A(C_4)$. The same arguments repeated for $u_2$ instead of $u_1$ prove this assertion.
\end{proof}

Being able to characterize all elements in an ear RAAG that arise from a degree $1$ or $2$ vertex in the underlying graph, we can rigorously define pearls.

\begin{definition}
Let $A(\Gamma)$ be an ear RAAG.
\begin{enumerate}
    \item An element $x \in A(\Gamma)$ is called a pearl if it satisfies the conditions of Lemma~\ref{lem:degtwovertices}, of Lemma~\ref{lem:degonevertices} or of Proposition~\ref{prop:specialcases}.
    \item It is called an end pearl if it arises from a degree $1$ vertex in a path (as specified in Lemma~\ref{lem:degonevertices} or in Proposition~\ref{prop:specialcases}) or if it is marked as a distinguished element in a cycle RAAG. All other pearls are called inner pearls.
    \item Two pearls $x$ and $y$ are called adjacent if they satisfy one of the following. They are inner pearls satisfying the conditions of Proposition~\ref{prop:adjacency}. Or $x$ is an end pearl and $y$ an inner pearl such that they fulfill the conditions of Lemma~\ref{lem:degonevertices} together with another inner vertex $z$. Or they satisfy one of the appropriate conditions of Proposition~\ref{prop:specialcases}.
\end{enumerate}
\end{definition}

\section{Vertex bases in ear RAAGs and loose ear decompositions}

First we characterize vertex bases in ear RAAGs without specifying that these elements are generators. This enables us then to characterize loose ear decompositions of graphs in terms of their associated RAAGs. The following proposition sets the base to assemble pearls into vertex bases.

\begin{proposition}\label{prop:basisextension}
Let $A(\Gamma)$ be an ear RAAG with $d(A(\Gamma)) \geq 5$ and $x, y, z \in A(\Gamma)$ be inner pearls such that $y, z \in C(x)$ and such that $y \langle x \rangle, z \langle x \rangle$ form a generating set of $C(x)/ \langle x \rangle$. Then $x$, $y$ and $z$ form part of a vertex basis such that the vertex corresponding to $x$ is adjacent to the one corresponding to $y$ as well as the one corresponding to $z$.
\end{proposition}

\begin{proof}
Let $x$ corresponds to a vertex $v_i$. According to the proof of Proposition~\ref{prop:adjacency}, there are words $a, b \in C(v_i)$ such that $y = av_j a^{-1}$ and $z = bv_k b^{-1}$ for $j, k \in \lbrace i-1, i, i+1 \rbrace$. Conjugating by $a^{-1}$ yields $y = v_j$, $x = v_i$ and $z = ev_k e^{-1}$ for $e \in C(v_i)$. We may assume that $e$ is a reduced word in $v_{i-1}^{\pm 1}, v_{i+1}^{\pm 1}$ that does not end in the letter $v_k$ or $v_k^{-1}$. Whether $y \langle x \rangle = v_j \langle v_i \rangle$ and $z \langle x \rangle = ev_k e^{-1} \langle v_i \rangle$ generate $\big\langle v_{i-1} \langle x \rangle, v_{i+1} \langle x \rangle \, \big\rangle \cong F_2$ is invariant under conjugation. Thus, we may further assume that the word $e$ does not start in $v_j$ or $v_j^{-1}$. If $e$ was still nontrivial, then $\lbrace y \langle x \rangle, z \langle x \rangle \rbrace$ would not a be generating set of $C(x)/ \langle x \rangle$. Therefore $e = 1$, $y = v_j$,  $z = v_k$ and $j \neq k$. We conclude that $\lbrace x, y, z \rbrace = \lbrace v_{i-1}, v_i, v_{i+1} \rbrace$.
\end{proof}

\begin{lemma}\label{lem:pathpearlchains}
Let $A(\Gamma)$ be a path RAAG. Let $\lbrace p_1, \dots, p_m \rbrace \subseteq A(\Gamma)$ be a set of pearls with the following properties. The triples $\lbrace p_1, p_2, p_3 \rbrace$ and $\lbrace p_{m-2}, p_{m-1}, p_m \rbrace$ satisfy the conditions of Lemma~\ref{lem:degonevertices} where $p_1$ and $p_m$ are the end pearls. The inner pearls $p_{i-1}$, $p_i$, $p_{i+1}$ satisfy the conditions of Proposition~\ref{prop:basisextension} for any $i \in \lbrace 3, \dots, m-2 \rbrace$. Then $\lbrace p_1, \dots, p_m \rbrace $ is a vertex basis.
\end{lemma}

\begin{proof}
If $v_1, \dots, v_n$ denote the vertices of $\Gamma$, then we may assume by Lemma~\ref{lem:degonevertices} that $p_1 = v_1$, $p_2 = v_2$ and $p_3 = v_3$. Since $p_3$ is an inner pearl, we see that $m, n \geq 4$. We remark that we have characterised vertex bases for $n \in \lbrace 2, 3 \rbrace$ in Proposition~\ref{prop:specialcases}. If $m \geq 5$, then $p_4$ is not an end pearl and we have at least one triple $\lbrace p_{i-1}$, $p_i$, $p_{i+1} \rbrace$ satisfying the conditions of Proposition~\ref{prop:basisextension} for $i \in \lbrace 3, \dots, m-2 \rbrace$. Assume by induction that $p_{i-1} = v_{i-1}$ and $p_i = v_i$. Following the proof of that proposition, there is $c \in \mathbb{Z}$ such that $p_{i+1} = v_{i-1}^c v_{i+1}v_{i-1}^{-c}$. By the treatment of RAAG automorphism in \cite[p.~498]{kob22}, we can apply the so-called partial conjugation that conjugates all vertex generators arising from the connected component of $A(\Gamma) \setminus \mathrm{star}(v_{i-1})$ containing $v_{i+1}$ by the element $v_{i-1}$ and that leaves all other vertex generators invariant. Since this completes the inductive step, we need to find an isomorphism of $A(\Gamma)$ that leaves $p_1, \dots, p_{m-1}$ invariant and maps $p_m$ to $v_m$. This is constructed in the proof of Lemma~\ref{lem:degonevertices}. Namely, in the first step, we deduce that $p_m = v_{m-2}^{\gamma}v_m^d v_{m-2}^{\delta}v_{m-1}^{\eta}$ for $\gamma, \delta \in \mathbb{Z}$ and $d = \pm 1$. The desired isomorphism is obtained by the dominated transvections and vertex inversions mapping $p_m$ to $v_m$.
\end{proof}

In Proposition~\ref{prop:specialcases} we have characterized vertex bases for cycle RAAGs $A(C_n)$ in the case where $n \in \lbrace 3, 4 \rbrace$.

\begin{lemma}\label{lem:cyclepearlchains}
Let $A(\Gamma)$ be a cycle RAAG with $d(A(\Gamma)) \geq 5$. Take $\lbrace p_1, \dots, p_m \rbrace \subseteq A(\Gamma)$ to be a set of pearls such that for any $i \in \lbrace 1, \dots, m \rbrace$ the triple $\lbrace p_{i-1}, p_i, p_{i+1} \rbrace$ satisfies the conditions of Proposition~\ref{prop:basisextension} where we consider the indices modulo $m$. Then $\lbrace p_1, \dots, p_m \rbrace$ is a vertex basis.
\end{lemma}

\begin{proof}
If $v_1, \dots, v_n$ denote the vertices of $\Gamma$, then we may assume that $p_1 = v_1$, $p_2 = v_2$ and $p_3 = v_3$ by Proposition~\ref{prop:basisextension}. Following the proof of this proposition with the triple $\lbrace p_2, p_3, p_4 \rbrace$, there is $c_2 \in \mathbb{Z}$ such that $p_4 = v_2^{c_2}v_4 v_2^{-c_2}$. Assume by induction that there are words of the form $w_i = \prod_{j = 2}^{i-2} v_j^{c_j}$ such that $p_i = w_i v_i w_i^{-1}$. Since $w_{i-1}v_{i-2}^{c_{i-2}} = w_i$ and $v_{i-2} \in C(v_{i-1})$, conjugating the triple $\lbrace p_{i-1}, p_i, p_{i+1} \rbrace$ by $w_i^{-1}$ yields $\lbrace v_{i-1}, v_i, w_i^{-1}p_{i+1}w_i \rbrace$. Using the latter triple, there is $c_{i-1} \in \mathbb{Z}$ such that $w_i^{-1}p_{i+1}w_i = v_{i-1}^{c_{i-1}}v_{i+1}v_{i-1}^{-c_{i-1}}$ by the proof of Proposition~\ref{prop:basisextension}. Conjugating the last equality by $w_i$ completes the inductive step. As $\Gamma$ is a cycle, we infer that
\[ \Big( \prod_{j = 2}^{m-1} v_j^{c_j}\Big) v_1 \Big( \prod_{j = 2}^{m-1} v_j^{c_j}\Big)^{-1} = p_1 = v_1 \, . \]
Conjugating by $v_2^{-c_2}$ and multiplying by $v_1^{-1}$, we deduce that
\[ \Big( \prod_{j = 3}^{m-1} v_j^{c_j}\Big) v_1 \Big( \prod_{j = 3}^{m-1} v_j^{c_j}\Big)^{-1}v_1^{-1} = 1 \, . \]
If $c_j \neq 0$ for any $j \in \lbrace 3, \dots, m-1 \rbrace$, then the left hand side would be a reduced word in $v_1^{\pm 1}, \dots, v_n^{\pm 1}$. As seen at the end of the proof of Proposition~\ref{prop:adjacency}, this word cannot equal the identity. This contradiction implies that $w_i = v_2^{c_2}$ and $p_i = v_2^{c_2}v_i v_2^{-c_2}$ for any $i \in \lbrace 4, \dots, m \rbrace$. Because $p_1, p_2, p_3 \in C(v_2)$, conjugating by $v_2^{-c_2}$ yields $p_i = v_i$ for every $i \in \lbrace 1, \dots, m \rbrace$. Finally, since the triple $\lbrace p_{m-1}, p_m, p_1 \rbrace$ satisfies the conditions of Proposition~\ref{prop:basisextension}, $v_m$ is adjacent to $v_1$ and $\lbrace p_1, \dots, p_m \rbrace$ is a vertex basis.
\end{proof}

Having been able to characterize vertex bases in ear RAAGs without specifying that the elements in question are vertex generators, we can make the following definition.

\begin{definition}
Let $A(\Gamma)$ be an ear RAAG. A pearl chain $\lbrace p_1, \dots, p_n \rbrace$ is a set of pearls in $A(\Gamma)$ that satisfies the conditions of Proposition~\ref{prop:specialcases}, of Lemma~\ref{lem:pathpearlchains} or of Lemma~\ref{lem:cyclepearlchains}. An ear RAAG together with a pearl chain $(A(\Gamma), \lbrace p_1, \dots, p_n \rbrace)$ is called a decorated ear RAAG.
\end{definition}

Pearl chains in turn enable us to characterize loose ear decompositions by RAAGs.

\begin{theorem}\label{thm:looseears}
Suppose $\big\lbrace (A(\Gamma_i), \lbrace p_{i, 1}, \dots, p_{i, n_i} \rbrace) \big\rbrace_{i = 1}^m$ is a finite sequence of decorated ear RAAGs. Define $G_1 := A(\Gamma_1)$ together with the distinguished elements $\lbrace p_{1,1}, \dots, p_{1, n_1} \rbrace$. Assuming that we have defined the group $G_k$ and have kept track of the distinguished elements $p_{i, j}$ for $i \in \lbrace 1, \dots, k \rbrace$ and $j \in \lbrace 1, \dots, n_i \rbrace$, we define $G_{k+1} := G_k \ast_{H_k} A(\Gamma_{k+1})$ where
\[ H_k = \lbrace 1 \rbrace, \, H_k = \lbrace p_{k+1, 1} = p_{i_1, j_1} \rbrace \text{ or } H_k = \lbrace p_{k+1, 1} = p_{i_1, j_1}, p_{k+1, n_{k+1}} = p_{i_2, j_2} \rbrace \, , \]
$i_1, i_2 \in \lbrace 1, \dots, k \rbrace$, $j_1 \in \lbrace 1, \dots, n_{i_1} \rbrace$ and $j_2 \in \lbrace 1, \dots, n_{i_2} \rbrace$ . Then $G_m$ is a RAAG whose underlying graph $\Gamma$ has the following loose ear decomposition. We define $\Gamma_1' := \Gamma_1$. If we denote by $v_{i,j}$ the vertex in $\Gamma_i$ corresponding to the pearl $p_{i, j} \in A(\Gamma_i)$, we construct the graph $\Gamma_{k+1}'$ from $\Gamma_k'$ and $\Gamma_{k+1}$ as a disjoint union, by identifying only $v_{k+1, 1}$ with $v_{i_1, j_1}$ or by identifying both $v_{k+1, 1}$ with $v_{i_1, j_1}$ and $v_{k+1, n_{k+1}}$ with $v_{i_2, j_2}$ depending on the corresponding move in the RAAGs. Finally, $\Gamma = \Gamma_m' = \bigcup_{i = 1}^m \Gamma_i$. \\

On the other hand, if we are given a loose ear decomposition of $\Gamma$ by the ears $\Gamma_1, \dots, \Gamma_m$, then the iterated amalgamated product of $A(\Gamma_1), \dots, A(\Gamma_m)$ constructed as above is isomorphic to the RAAG $A(\Gamma)$.
\end{theorem}

\begin{proof}
Recall that pearl chains are vertex bases of ear RAAGs. Thus, in the above iterated amalgamated products, we are identifying vertex generators of two RAAGs, giving the presentation of another RAAG. The resulting group presentation is in bijective correspondence with the associated loose ear decomposition.
\end{proof}

\begin{remark}
Ear decompositions are the special case of loose ear decompositions where the first ear is a cycle and the endpoints of each subsequent ear are contained the previous ears. Hence, we can reformulate the above theorem for ear decompositions if we use the appropriate identifications and amalgamated products.
\end{remark}

\begin{definition}
A (loose) ear decomposition of a RAAG $A(\Gamma)$ is a specification of ear RAAGs $A(\Gamma_1)$, \dots, $A(\Gamma_m)$ and of an iterated amalgamated product representing $A(\Gamma)$ as in Theorem~\ref{thm:looseears}. We call the set of elements $\lbrace p_{i, j} \mid 1 \leq i \leq m, 1 \leq j \leq n_i \rbrace$ of $A(\Gamma)$ the pearls of the (loose) ear decomposition.
\end{definition}

In the same way as there are many possible choices of vertex bases for a RAAG, there is no unique choice of a pearl chain for an ear RAAG. In particular, (loose) ear decomposition of RAAGs are not unique. However, it is unlikely that there are algorithms using loose ear decompositions that determine for a RAAG whether its underlying graph possesses any graph theoretic property of interest. Namely, unless the RAAG is given as an iterated amalgamated product of the desired form, one needs to determine its underlying graph in order to find a loose ear decomposition of the RAAG.

\section{Applications of loose ear decompositions of RAAGs}

We provide various applications of loose ear decompositions of RAAGs in this section. Our first application is a characterization of graph subdivisions in terms of RAAGs. For this, it is more convenient to consider the reverse procedure of an edge subdivision. Namely, a vertex smoothing is the process of removing a degree 2 vertex from a graph and connecting its former neighbors by an edge~\cite[p.~305]{gro19}.

\begin{proposition}
Suppose that $A(\Gamma)$ is a RAAG. If we are given a loose ear decomposition of $A(\Gamma)$, then a vertex smoothing of $v \in V(\Gamma)$ determines a homomorphism of RAAGs $A(\Gamma) \rightarrow A(\Gamma)/ \langle\! \langle p, [q, r] \rangle\! \rangle$ where $p$ is the pearl corresponding to $v$ and $q$, $r$ denote its adjacent pearls. On the other hand, let $p'$ be a pearl in a loose ear decomposition of $A(\Gamma)$ and $q'$, $r'$ be its unique adjacent pearls. Then the homomorphism $A(\Gamma) \rightarrow A(\Gamma)/ \langle\! \langle p', [q', r'] \rangle\! \rangle$ determines a vertex smoothing of $\Gamma$.
\end{proposition}

\begin{proof}
Given a loose ear decomposition of $A(\Gamma)$, denote by $p$ the pearl corresponding to $v$ and by $q$ and $r$ its adjacent pearls. Because the RAAG $A(\Gamma)$ can be defined by a group presentation using the pearls of the loose ear decomposition, smoothing out the vertex $v$ corresponds to the homomorphism of RAAGs $A(\Gamma) \rightarrow A(\Gamma)/ \langle\! \langle p, [q, r] \rangle\! \rangle$. For the other implication, the vertex $v'$ associated to the pearl $p'$ has two unique neighbors. The homomorphism $A(\Gamma) \rightarrow A(\Gamma)/ \langle\! \langle p', [q', r'] \rangle\! \rangle$ corresponds to smoothing out the vertex $v'$.
\end{proof}

\begin{definition}
A homomorphism $\varphi: A(\Gamma) \rightarrow A(\Lambda)$ as in the previous proposition is called a smoothing homomorphism.
\end{definition}

\begin{lemma}\label{lem:subdivisions}
Let $A(\Gamma)$ and $A(\Lambda)$ be RAAGs. Then $\Gamma$ is a graph subdivision of $\Lambda$ if and only if $A(\Gamma) \cong A(\Lambda)$ or there exists a finite sequence of smoothing homomorphisms $\lbrace \varphi_{i}: A(\Gamma_{i}) \rightarrow A(\Gamma_{i-1}) \rbrace_{i = 1}^m$ such that $A(\Gamma_m) = A(\Gamma)$ and $A(\Gamma_0) \cong A(\Lambda)$.
\end{lemma}

\begin{proof}
The case where the graphs $\Gamma$ and $\Lambda$ are isomorphic corresponds to the case where $A(\Gamma)$ and $A(\Lambda)$ are isomorphic, which is covered by~\cite[Theorem~15.2.6]{kob22}. Assume that the graphs $\Gamma$ and $\Lambda$ are not isomorphic. Every smoothing homomorphism corresponds to a vertex smoothing which corresponds uniquely to an edge subdivision. Hence, for a graph subdivision of $\Gamma$, we can construct the corresponding sequence of smoothing homomorphisms. If we are given a sequence of smoothing homomorphisms instead, we invoke \cite[Theorem~15.2.6]{kob22} to obtain the underlying edge subdivisions forming the graph subdivision.
\end{proof}

Graph subdivisions can be used to characterize planar graphs in terms of RAAGs. For this, let $K_5$ be the complete graph on $5$ vertices and $K_{3, 3}$ be the complete bipartite graph obtained as a join of two graphs of $3$ isolated vertices each.

\begin{theorem}\label{thm:planarity}
Let $A(\Gamma)$ be a RAAG. Assume that $\lbrace \varphi_{i}: A(\Gamma_{i}) \rightarrow A(\Gamma_{i-1}) \rbrace_{i = 1}^m$ is a sequence of smoothing homomorphisms such that $A(\Gamma_m) = A(\Gamma)$ and such that there is no smoothing homomorphism from $A(\Gamma_0)$ to any other RAAG. Then the graph $\Gamma$ is planar if and only $A(\Gamma_0)$ does not contain $\mathbb{Z}^5$ as a subgroup and does not possess a loose ear decomposition with pearls $\lbrace p_i, q_i \rbrace_{i = 1}^3$ such that the subgroups $P := \langle p_1, p_2, p_3 \rangle$ and $Q := \langle q_1, q_2, q_3 \rangle$ are isomorphic to $F_3$ and $\langle P, Q \rangle \cong F_3 \oplus F_3$.
\end{theorem}

\begin{proof}
As $\Gamma$ is a subdivision of $\Gamma_0$, the former graph is planar if and only if the latter one is~\cite[Proposition~7.5.15]{gro19}. The graph $\Gamma_0$ is not a subdivision of another graph by our assumptions. This together with Kuratowski's Theorem \cite[Theorem~7.4.1]{gro19} implies that $\Gamma_0$ is planar if and only if it does not contain a copy of $K_5$ or $K_{3, 3}$ as a subgraph. According to \cite[Theorem~15.3.5]{kob22}, $\Gamma_0$ contains $K_5$ as a subgraph if and only if $A(\Gamma_0)$ contains $\mathbb{Z}^5$ as a subgroup. If $u_i$ are the vertices corresponding to the pearls $p_i$ and $v_i$ the ones corresponding to $q_i$, then $P$ (respectively $Q$) is isomorphic to $F_3$ if and only if there are no edges between the vertices $u_i$ (respectively $v_i$) by \cite[Theorem~15.3.2]{kob22}. Using \cite[Theorem~15.3.1]{kob22}, $\langle P, Q \rangle$ decomposes as non-trivial direct product only if the vertices $u_i$ form a non-trivial join with the vertices $v_j$, meaning a copy of $K_{3,3}$.
\end{proof}

Moving forward, the following equivalent formulation of graph minors allows us to describe them in terms of RAAGs. An edge contraction of a graph consists of identifying two adjacent vertices while deleting the edge between them~\cite[p.~21]{die17}.

\begin{proposition}
Let $X$ and $Y$ be graphs. Then $X$ is a minor of $Y$ if and only if there are finite sequences of graphs $\lbrace A_i \rbrace_{i = 1}^a$, $\lbrace B_j \rbrace_{j = 1}^b$ and $\lbrace C_k \rbrace_{k = 1}^c$ with $A_1 = Y$, $A_a = B_1$, $B_b = C_c$, $C_1 = X$ and the following properties. The graph $A_{i+1}$ is obtained from $A_i$ by a vertex deletion, $B_{j+1}$ from $B_j$ by an edge contraction and $C_k$ from $C_{k+1}$ by an edge deletion.
\end{proposition}

We can visualize the sequences of graphs as
\begin{center}
\begin{tikzcd}
{} & & & {} & & & X = C_1 \\
& & & & & & \\
Y = A_1 \arrow[rrr, rightsquigarrow, "\text{vertex deletions}"] & & & A_a = B_1 \arrow[rrr, rightsquigarrow, "\text{edge contractions}"] & & & B_b = C_c \arrow[uu, rightsquigarrow, "\text{edge deletions}"]
\end{tikzcd}
\end{center}

\begin{proof}
It follows from~\cite[Corollary~1.7.2]{die17} that $\lbrace A_i \rbrace_{i = 1}^a$, $\lbrace B_j \rbrace_{j = 1}^b$ and $\lbrace C_k \rbrace_{k = 1}^c$ give rise to a graph minor. On the other hand, assume that $f: S \subseteq V(Y) \rightarrow V(X)$ is a function describing $X$ as a minor of $Y$. We can retrieve the subset $S$ of $V(Y)$ by a sequence of vertex deletions. For any $x \in X$, edge contractions ensure to contract every subgraph $f^{-1}(x)$ to a single vertex. We then obtain $X$ by a sequence of edge deletions.
\end{proof}

If we choose a loose ear decomposition, the vertex deletion $A_i \rightsquigarrow A_{i+1}$ corresponds to a homomorphism $A(A_i) \rightarrow A(A_i)/ \langle\! \langle p \rangle\! \rangle$ where $p$ is a pearl and $A(A_i)/ \langle\! \langle p \rangle\! \rangle \cong A(A_{i+1})$. We call such a homomorphism a pearl deletion. The edge contraction $B_j \rightsquigarrow B_{j+1}$ corresponds to a homomorphism $A(B_j) \rightarrow A(B_j)/ \langle\! \langle pq^{-1} \rangle\! \rangle$ where $p$ and $q$ are adjacent pearls and $A(B_j)/ \langle\! \langle pq^{-1} \rangle\! \rangle \cong A(B_{j+1})$. As we identify $p$ and $q$, we call such a homomorphism an identification homomorphism. Lastly, the edge deletion $C_{k+1} \rightsquigarrow C_k$ corresponds to a homomorphism $A(C_k) \rightarrow A(C_k)/ \langle\! \langle [p, q] \rangle\! \rangle$ where $p$ and $q$ are not adjacent and $A(C_k)/ \langle\! \langle [p, q] \rangle\! \rangle \cong A(C_{k+1})$. Since we impose that two elements commute, we call such a homomorphism an abelianizing homomorphism. We reformulate the above proposition as

\begin{theorem}\label{thm:minors}
Let $A(\Gamma)$ and $A(\Lambda)$ be RAAGs. Then $\Lambda$ is a minor of $\Gamma$ if and only if there exists a sequence of pearl deletions $\lbrace A(A_i) \rightarrow A(A_{i+1}) \rbrace_{i = 1}^a$, one of identification homomorphisms $\lbrace A(B_j) \rightarrow A(B_{j+1}) \rbrace_{j = 1}^b$ and one of abelianizing homomorphisms $\lbrace A(C_k) \rightarrow A(C_{k+1}) \rbrace_{k = 1}^c$ such that $A(\Gamma) \cong A(A_1)$, $A(A_a) \cong A(B_1)$, $A(B_b) \cong A(C_c)$ and $A(C_1) \cong A(\Lambda)$.
\end{theorem}

This can be pictorially summarized as
\begin{center}
\begin{tikzcd}[column sep = scriptsize]
{} & & & {} & & & & A(\Lambda) = A(C_1) \arrow[dd, "\text{abelianizing homom's}" description] \\
& & & & & & & \\
A(\Gamma) = A(A_1) \arrow[rrr, "\text{pearl deletions}"] & & & A(A_a) = A(B_1) \arrow[rrrr,"\text{identifying homom's}"] & & & & A(B_b) = A(C_c)
\end{tikzcd}
\end{center}

\begin{definition}
If two RAAGs $A(\Gamma)$ and $A(\Lambda)$ satisfy the conditions of the previous theorem, then $A(\Lambda)$ is a called a RAAG minor of $A(\Gamma)$.
\end{definition}

The key of characterizing minor-closed properties in terms of RAAGs is the following theorem corresponding to Theorem~12.7.1 and Corollary~12.7.2 in~\cite{die17}:

\begin{theorem}[Robertson-Seymour]
Any minor-closed graph property can be characterized by a finite set of forbidden minors.
\end{theorem}

All our applications of RAAGs to graph minors are based on this theorem. For instance,

\begin{proposition}\label{prop:forests}
Let $A(\Gamma)$ be a RAAG. Then $\Gamma$ is a forest if and only if $A(\Gamma)$ does not contain $\mathbb{Z}^3$ as a RAAG minor.
\end{proposition}

\begin{proof}
A graph is a forest if it does not contain a cycle, meaning a subdivision of $C_3$ where $A(C_3) = \mathbb{Z}^3$.
\end{proof}

\begin{theorem}
Let $A(\Gamma)$ be a RAAG. Then $\Gamma$ is planar if and only if $A(\Gamma)$ does not contain $\mathbb{Z}^5$ or $F_3 \oplus F_3$ as RAAG minors.
\end{theorem}

\begin{proof}
According to Wagner's Theorem~\cite[Theorem~4.4.6]{die17}, a graph is planar if and only it does not contain $K_5$ or $K_{3, 3}$ as a graph minor where $A(K_5) =  \mathbb{Z}^5$ and $A(K_{3, 3}) = F_3 \oplus F_3$.
\end{proof}

For any embedding of a planar graph $\Gamma$ into $\mathbb{R}^2$ without edge crossings, the connected components of $\mathbb{R} \setminus \Gamma$ are called faces~\cite[p.~92]{die17}. The graph $\Gamma$ is called outerplanar if all its vertices lie on the boundary of the unique unbounded face~\cite[p.~115]{die17}.

\begin{lemma}\label{lem:outerplanar}
Let $A(\Gamma)$ be a RAAG. Then $\Gamma$ is outerplanar if and only if it does not contain $\mathbb{Z}^4$ or $F_2 \oplus F_3$ as RAAG minors.
\end{lemma}

\begin{proof}
A graph $\Gamma$ is outerplanar if and only if it does not contain the complete graph $K_4$ or the complete bipartite graph $K_{2, 3}$ as graph minors~\cite[p.~115]{die17}, where $A(K_4) = \mathbb{Z}^4$ and $A(K_{2, 3}) = F_2 \oplus F_3$.
\end{proof}

In this manner, one could characterize more graph theoretic properties in terms of RAAGs. For example, a graph is embeddable in a given surface $S$ if and only if it does not contain any of a finite set of forbidden minors~\cite[Corollary~12.7.3]{die17}. As soon as such a finite set is known, one can write down the associated RAAGs. Furthermore, a graph $\Gamma$ is called an apex graph if there exists a vertex $v$ such that $\Gamma \setminus \lbrace v \rbrace$ is planar~\cite[p.~809]{gup91}. It is known that apex graphs can be also characterized by a finite set of forbidden minors, but it is unknown what these are~\cite[p.~809]{gup91}. As soon as this is known, one could characterize apex graphs by looking at their RAAGs. \\

In our last applications, we specifically want to restrict our attention to ear decompositions. Let us introduce two notions of connectivity for graphs. A graph is called $k$-vertex-connected if removing at most $k-1$ vertices does not disconnect it, while removing $k$ (specific) vertices does~\cite[p.~237]{schr04}. Similarly, a graph is $k$-edge-connected if removing at most $k-1$ edges does not disconnect it, but removing $k$ (particular) edges does~\cite[p.~237]{schr04}.

\begin{theorem}\label{thm:twoedgeconn}
Let $\Gamma$ be a graph. Then $\Gamma$ is $2$-edge-connected if and only if the RAAG $A(\Gamma)$ admits an ear decomposition.
\end{theorem}

\begin{proof}
This follows from Theorem~15.17 in~\cite{schr04}.
\end{proof}

According to page~252 of~\cite{schr04}, we call an ear decomposition $\lbrace \Gamma_i \rbrace_{i = 1}^m$ of a graph $\Gamma$ proper if for any $i > 1$ the ear $\Gamma_i$ is a path. Call the corresponding decomposition of $A(\Gamma)$ also a proper ear decomposition.

\begin{lemma}\label{lem:twovertexconn}
Let $\Gamma$ be a graph such that $d(A(\Gamma)) \geq 2$. Then $\Gamma$ is $2$-vertex-connected if and only if $A(\Gamma)$ has a proper ear decomposition.
\end{lemma}

\begin{proof}
Since $|V(\Gamma)| = d(A(\Gamma)) \geq 2$ by Proposition~\ref{prop:degsandvertices}, this follows from Theorem~15.16 in~\cite{schr04}.
\end{proof}

After these connectivity results, let us consider series-parallel graphs that can be defined as follows according to~\cite[p.~41]{epp92}. A graph $\Gamma$ consisting of a single edge is a series-parallel graph where we call its two vertices $s_{\Gamma}$, $t_{\Gamma}$ its terminals. If $(X, s_X, t_X)$ and $(Y, s_Y, t_Y)$ are two series parallel graphs, then we can form a new series parallel graph $(\Gamma, s_{\Gamma}, t_{\Gamma})$ by identifying $s_X$ with $s_Y$, declaring this to be $s_{\Gamma}$ and doing the same for $t_X$, $t_Y$ and $t_{\Gamma}$. We could equally define a different series-parallel graph $(\Lambda, s_{\Lambda}, t_{\Lambda})$ by identifying $t_X$ with $s_Y$ and declaring the terminals to be $s_X$ and $t_Y$. These graphs are related to a particular form of ear decomposition. A nested ear decomposition of a graph $\Gamma$ is a loose ear decomposition $\lbrace \Gamma_i \rbrace_{i = 1}^m$ consisting of paths and satisfying the following conditions. For any $1 < i$ there is $1 \leq j < i$ such that the endpoints of $\Gamma_i$ are contained in $\Gamma_j$. If  the endpoints of $\Gamma_i$ and $\Gamma_j$ lie in $\Gamma_k$, then the subpath in $\Gamma_k$ between the endpoints of $\Gamma_i$ contains either both or none of endpoints of $\Gamma_j$~\cite[p.~43]{epp92}. Contrary to ear decompositions, the first ear of a nested ear decomposition is not a cycle. Call the corresponding loose ear decomposition of $A(\Gamma)$ a nested ear decomposition. We reformulate Theorem~1 in~\cite{epp92} as

\begin{lemma}\label{lem:seriesparallel}
Let $\Gamma$ be a RAAG. Then $\Gamma$ is a series-parallel graph if and only if $A(\Gamma)$ has a nested ear decomposition. In this case, the terminals of $\Gamma$ correspond to the end pearls of the first ear RAAG in the nested ear decomposition of $A(\Gamma)$.
\end{lemma}

A graph $\Gamma$ is called factor-critical if for every vertex $v \in V(\Gamma)$ the graph $\Gamma \setminus \lbrace v \rbrace$ has a perfect matching~\cite[p.~425]{schr04}, meaning a subset $M \subseteq E(\Gamma \setminus \lbrace v \rbrace)$ such that any two elements of $M$ are disjoint and every vertex in $\Gamma \setminus \lbrace v \rbrace$ is covered by at least an edge in $M$~\cite[p.~23]{schr04}. Then Theorem~24.9 in~\cite{schr04} can be rephrased as

\begin{lemma}\label{lem:factorcritical}
Let $\Gamma$ be a graph. Then $\Gamma$ is factor-critical if and only if $A(\Gamma)$ has an ear decomposition $\lbrace A(\Gamma_i) \rbrace_{i = 1}^m$ such that $d(A(\Gamma_i))$ is odd if $A(\Gamma_i)$ is a cycle RAAG and even if $A(\Gamma_i)$ is a path RAAG.
\end{lemma}

Lastly, for an ear RAAG $A(\Gamma)$, set the quantity $E(A(\Gamma)) = \lfloor \frac{1}{2} d(A(\Gamma)) \rfloor$ if $\Gamma$ is a cycle and $E(A(\Gamma)) = \lfloor \frac{1}{2} (d(A(\Gamma))-1) \rfloor$ if $\Gamma$ is a path. One can think of $E(A(\Gamma))$ as roughly half the edges of $\Gamma$.

\begin{proposition}\label{prop:maxjoinsize}
Let $\Gamma$ be a $2$-vertex-connected graph. Then the following hold.
\begin{enumerate}
    \item If $\lbrace A(\Gamma_i) \rbrace_{i = 1}^m$ is any ear decomposition of $A(\Gamma)$, then $m = |E(\Gamma)|-|V(\Gamma)|+1$. In terms of group cohomology over a field $F$, this is equivalent to
    \[ m = \mathrm{dim}_F \big(H_F^2(A(\Gamma), F)\big)- \mathrm{dim}_F \big(H_F^1(A(\Gamma), F)\big)+1 \, . \]
    \item The maximum size of a join in $\Gamma$ equals
    \[ \mathrm{min} \Bigg\lbrace \sum_{i = 1}^m E(A(\Gamma_i)) \mid \lbrace A(\Gamma_1) \rbrace_{i = 1}^m \text{ ear decomposition of } A(\Gamma) \Bigg\rbrace \, . \]
\end{enumerate}
\end{proposition}

\begin{proof}
The first part is due to~\cite[p.~511]{schr04} and \cite[Proposition~15.2.4]{kob22}, where the second part is a reformulation of~\cite[Theorem~29.11]{schr04}.
\end{proof}

\section*{Acknowledgement}

I am indebted to the referee whose invaluable feedback has led to an improvement of this paper. I am also thankful to Thomas Koberda for pointing out that it is unlikely that the results of this paper could be algorithmically implemented. Lastly, I would like to thank Ashot Minasyan as well as Jens M.\! Schmidt for a useful comment.

\addcontentsline{toc}{section}{References}
\renewcommand{\bibname}{References}
\bibliographystyle{plain}  
\bibliography{refs}        

Department of Mathematics, University of British Columbia, 125-1984 Mathematics Road, Vancouver, BC, V6T 1Z2, Canada \newline
\textit{Email address}: gheorghiu.max@gmail.com \newline
\textit{URL}: https://www.math.ubc.ca/user/2675

\end{document}